\documentclass[twoside,11pt]{amsart}

\usepackage{amssymb}
\usepackage[leqno]{amsmath}
\usepackage{amsfonts}
\usepackage{amsopn}
\usepackage{amstext}
\usepackage{amsthm}

\providecommand{\cal}{\mathcal}
\renewcommand{\Bbb}{\mathbb}

\newenvironment{pf}{\begin{proof}}{\end{proof}}

\newcommand{\Cee}{{\cal{C}}}

\newcommand{\Tau}{{\cal{T}}}

\newcommand{\Err}{{\Bbb{R}}}

\newcommand{\lam}{{\lambda}}
\newcommand{\al}{\alpha}

\newcommand{\sig}{\sigma}
\newcommand{\eps}{\varepsilon}
\renewcommand{\phi}{\varphi}
\renewcommand{\rho}{\varrho}

\newcommand{\rest}{\restriction}

\newcommand{\ntr}{n\in\omega}

\newcommand{\loe}{\leqslant}
\newcommand{\goe}{\geqslant}

\newcommand{\subs}{\subseteq}
\newcommand{\sups}{\supseteq}
\newcommand{\nnempty}{\ne\emptyset}

\newcommand{\dist}{\operatorname{dist}}

\newcommand{\conv}{\operatorname{conv}}

\newcommand{\dom}{\operatorname{dom}}
\newcommand{\rng}{\operatorname{rng}}

\newcommand{\poset}{{\Bbb{P}}}

\newcommand{\Es}{{\cal{S}}}

\newcommand{\Land}{\;\&\;}

\newcommand{\Lev}{\operatorname{Lev}}
\newcommand{\Ht}{\operatorname{ht}}

\newcommand{\length}{\operatorname{length}}
\newcommand{\concat}{{}^\smallfrown}

\newcommand{\by}{/_}

\newtheorem{tw}{Theorem}[section]
\newtheorem{wn}[tw]{Corollary}
\newtheorem{lm}[tw]{Lemma}
\newtheorem{prop}[tw]{Proposition}
\newtheorem{claim}[tw]{Claim}
\theoremstyle{definition}
\newtheorem{df}[tw]{Definition}

\theoremstyle{remark}
\newtheorem{uwgi}[tw]{Remark}

\newcommand{\setof}[2]{\{#1\colon #2\}}

\newcommand{\sn}[1]{\{#1\}} 
\newcommand{\dn}[2]{\{#1,#2\}} 
\newcommand{\map}[3]{#1\colon #2 \to #3} 
\newcommand{\img}[2]{#1[#2]} 

\newcommand{\baire}{\omega^\omega}
\newcommand{\app}{\operatorname{app}}

\newcommand{\rk}{\operatorname{rk}}
\newcommand{\suppt}{\operatorname{suppt}}

\newcommand{\biglam}{\lam_{\omega_1}(\aleph_0)}
\newcommand{\colors}{\operatorname{colors}}

\newcommand{\jezyk}{{\mathcal L}_{\omega,\omega_1}}
\newcommand{\Dee}{{\mathcal D}}

\title{Analytic Colorings}

\author{Wies{\l}aw Kubi\'s}
\address{Department of Mathematics,
Ben-Gurion University of the Negev, Beer-Sheva, Israel
\\ \textit{and} \\ Institute of Mathematics, University of Silesia,
Katowice, Poland}
\email{kubis@math.bgu.ac.il}

\author{Saharon Shelah $^{(\dagger)}$}
\address{Institute of Mathematics, Hebrew University of
        Jerusalem, Israel \\ \textit{and} \\ Department of Mathematics Rutgers
        University, New-Brunswick} \email{shelah@math.huji.ac.il}

\thanks{$(\dagger)$ This research was supported by The Israel Science Fundation.
Publication 802.}

\begin{document}


\begin{abstract}
We investigate the existence of perfect homogeneous sets for analytic colorings. An {\em analytic coloring} of $X$ is an analytic subset of $[X]^N$, where $N>1$ is a natural number. We define an absolute rank function on trees representing analytic colorings, which gives an upper bound for possible cardinalities of homogeneous sets and which decides whether there exists a perfect homogeneous set. We construct universal $\sig$-compact colorings of any prescribed rank $\gamma<\omega_1$.
These colorings consistently contain homogeneous sets of cardinality $\aleph_\gamma$ but they do not contain perfect homogeneous sets. As an application, we discuss the so-called defectedness coloring of subsets of Polish linear spaces.

{\bf 2000 Mathematics Subject Classification:} Primary: 03E05, 03E15. Secondary: 03E35, 54H05.

{\bf Keywords and phrases:} Analytic coloring, tree, homogeneous set, rank of a coloring tree.
\end{abstract}

\maketitle

\section{Introduction}

A classical result, due to Galvin \cite{Galvin} (see also
\cite[p. 130]{Kechris}), says that if $X$ is a nonempty Polish space,
$[X]^2=C_0\cup C_1$, where both $C_0$, $C_1$ are Borel, then there exists a
perfect set $P\subs X$ which is homogeneous for $C_i$, i.e. $[P]^2\subs
C_i$, for some $i<2$. As it was noticed by Blass in \cite{BlassEx}, this
result is no longer valid for colorings of triples. It is natural to ask
for sufficient conditions under which a given Borel (or, more generally,
analytic) coloring $C\subs [X]^N$ has a perfect homogeneous set (i.e. a
perfect set $P$ such that $[P]^N\subs C$). 

It is well known that an analytic
set is either countable or contains a perfect set. Another result of this nature is Silver's theorem: a co-analytic equivalence relation has either countably many equivalence classes or has a perfect set of pairwise unrelated elements.
A natural conjecture
could be as follows: an analytic coloring which contains an uncountable
homogeneous set, contains also a perfect one. This conjecture is true for
$G_\delta$ colorings of Polish spaces (see \cite{Sh522} or \cite{Kubis1}),
but in general it is false. There exist (in ZFC) $\sig$-compact colorings
which contain uncountable but not perfect homogeneous sets. However, it is
consistent that continuum is arbitrarily large and every analytic coloring
which contains a homogeneous set of size $\aleph_2$, contains also a perfect
one (see \cite{Sh522}). This topic has been investigated by the second
author in \cite{Sh522}, where an increasing continuous sequence of cardinals
$\setof{\lam_\al(\aleph_0)}{\al\loe\omega_1}$ is defined, using a certain
model-theoretic rank, and it is shown that every analytic coloring with a
homogeneous set of size $\biglam$ has also a perfect homogeneous set. Examples of $\sig$-compact colorings with large but not perfect homogeneous sets are constructed by means of a ccc forcing. For more detailed statements of these results, see Section \ref{krowka}.

In this paper we investigate analytic colorings without perfect homogeneous sets. We define a rank function on trees inducing analytic colorings. The rank of such a tree is either a countable ordinal or $\infty$, in the latter case the coloring has a perfect homogeneous set. We show that if an analytic coloring has a homogeneous set of size $\lam_\gamma(\aleph_0)$ then the rank of any tree representing this coloring is $\goe\gamma$. For a given ordinal $\gamma<\omega_1$ we construct a ``universal" $\sig$-compact pair coloring $C_\gamma$ which is represented by a tree of rank $\gamma$ and every $\sig$-compact coloring defined by a suitable tree of rank $\loe\gamma$ can be embedded into $C_\gamma$. If $\gamma>0$ then $C_\gamma$ contains an uncountable homogeneous set. Furthermore, for any cardinal $\lam\loe\lam_\gamma(\aleph_0)$ we construct a ccc forcing $\poset$ of size $\lam$ which adds a $C_\gamma$-homogeneous set of cardinality $\lam$. In fact, such a set is defined from a $\Dee$-generic filter, where $\Dee$ is some family of $\lam$ many dense subsets of $\poset$. Thus if $MA_{\lam_\gamma(\aleph_0)}$ holds then $C_\gamma$ contains a homogeneous set of cardinality $\lam_\gamma(\aleph_0)$.
The coloring $C_\gamma$ is constructed in ZFC, in fact it is induced by a certain tree which is built by using finite objects. This improves a result from \cite{Sh522}, where a ccc forcing notion is constructed, which adds generically a $\sig$-compact coloring with a perfect homogeneous set of size $\lam_\gamma(\aleph_0)$ but without a perfect one.

Finally, as an application we discuss the {\em defectedness coloring} of a subset $X$ of a Polish linear space: $C=\setof{S\in[X]^N}{\conv S\not\subs X}$. A homogeneous set for this coloring is called a {\em clique in} $X$. Defectedness colorings were studied in \cite{FKo, GKKS, Ko1, KPS, Kubis1}. We show that the defectedness coloring of a Borel set is analytic and we prove that every $F_\sig$ coloring of the Cantor set can be ``realized" as defectedness of some bounded $G_\delta$ set in $\Err^{N+1}$.

The paper is organized as follows. Section \ref{pralinki} contains definitions and notation, including the model-theoretic rank from \cite{Sh522}, which is used to define cardinals $\lam_\gamma(\aleph_0)$. Section \ref{krowka} contains basic properties of analytic colorings and the statements of known results. In Section \ref{rangakoloru} we define the rank of a tree representing an analytic coloring; we state some of its basic properties and we show the relation to cardinals $\lam_\gamma(\aleph_0)$. In Section \ref{uniwers} we construct the mentioned universal $\sig$-compact colorings and the forcing notions adding homogeneous sets.
In the last section we discuss defectedness colorings of Borel subsets of Polish linear spaces.

\section{Preliminaries}\label{pralinki}

In this section we give the necessary definitions and notation used in the paper.

\subsection{Notation}
According to \cite{Jech}, a tree is a partially ordered set $(T,\loe)$ such that for every $x\in T$ the set $\setof{y\in T}{y<x}$ is well ordered. However, we will consider only trees of height $\loe\omega$, so for us a {\em tree} is a partially ordered set $(T,\loe)$ such that $\setof{y\in T}{y<x}$ is finite for every $x\in T$. We will use standard notation: $\Lev_n(T)$ is the $n$-th level of $T$. A {\em branch through} $T$ is an infinite chain in $T$. We denote by $[T]$ the set of all branches through $T$. A tree $T$ is {\em ever-branching} if for every $k<n<\Ht(T)$, for every $x\in\Lev_k(T)$ there is $y\in \Lev_n(T)$ with $x<y$. We will usually consider subtrees of $\omega^{<\omega}$ or, more generally, $\bigcup_{\ntr}(\omega^n)^k$. Every closed subset $A$ of the Baire space $\baire$ is of the form $[T]$ where $T=\setof{x\rest n}{x\in A,\;\ntr}$.

A metric space $X$ is {\em analytic} if it is a continuous image of the Baire space $\baire$. If $X$ is a subset of a Polish space $Y$ then $X$ is analytic iff $X=\setof{y\in Y}{(\exists\;t\in\baire)\;(y,t)\in F}$ for some closed set $F\subs Y\times \baire$. If $Y=\baire$ then $F=[T]$ for some tree $T$ consisting of pairs $(s,t)$, where $s,t\in\omega^n$, $\ntr$. We will say that $X$ is {\em induced} by $T$. 

Let $N>1$ be a natural number. An {\em $N$-coloring} of a space $X$ is a subset of $[X]^N$. Denoting by $X^{(N)}$ the subspace of $X^N$ consisting of all one-to-one $N$-tuples, we identify $[X]^N$ with the quotient $X^{(N)}\by\sim$, where $\sim$ is the obvious equivalence relation (which identifies $N$-tuples with the same range).
Thus we may speak about analytic, Borel, closed colorings, etc. A set $K$ is {\em $C$-homogeneous} (or {\em homogeneous for} $C$) if $[K]^N\subs C$. Write $\suppt(C)=\bigcup C$. Observe that, looking at $C$ as a subset of $X^N$, $\suppt(C)$ is the projection of $C$ onto $X$ (fixing one arbitrary coordinate). Thus, if $C$ is analytic then so is $\suppt(C)$.

\subsection{Coloring trees}

Let $C\subs [\baire]^N$ be an analytic $N$-coloring. Then, for some closed set $B\subs[\baire]^N\times\baire$ we have $C=\setof{v\in[\baire]^N}{(\exists\;\tau)\;(v,\tau)\in B}$. Let $T=\setof{(v\rest n,t\rest n)}{(v,t)\in B}$, where $v\rest n=\setof{x\rest n}{x\in v}$. Then $T$ is a tree such that $[T]=B$. Every analytic coloring is induced by a tree of this form.
Below we give a formal definition of trees which represent analytic colorings of $\baire$. 

\begin{df} Fix $N<\omega$. We define a partial order $<$ on the set $B=\bigcup_{\ntr}[\omega^n]^N\times\omega^n$ by setting $(v,t)<(v',t')$ iff $t\subs t'$ and $v=v'\rest n$ where $n=\dom(t)$. Given $T\subs B$ we define the {\em $n$-th level of\/ $T$}, $\Lev_n(T)=T\cap ([\omega^n]^N\times \omega^n)$ and we write $n=\length(v,t)$ if $(v,t)\in[\omega^n]^N\times \omega^n$. Define $\Ht(T)=\min\setof{n\loe\omega}{(\forall\;m\goe n)\;\Lev_m(T)=\emptyset}$ (so $\Ht(T)\loe\omega$).

An {\em $N$-coloring tree} is a set $T\subs B$ satisfying the following condition: 
\begin{enumerate}
	\item[($\bullet$)] If $(v,t)\in \Lev_n(T)$ and $n<m<\Ht(T)$ then there exists $(v',t')\in \Lev_m(T)$ such that $(v,t)<(v',t')$.
\end{enumerate}
Note that we do not require $(v\rest n, t\rest n)\in T$ for $n<\length(v,t)$, even if $|v\rest n|=|v|$. A coloring tree is a tree in the usual sense, but the notion of level is not the same as in the general sense of the tree -- a coloring tree $T$ may consists of elements of length $\goe n$ and in this case $\Lev_n(T)$ is in fact the least level of $T$. 
A coloring tree $T$ of infinite height induces an analytic coloring of the Baire space:
$$C=\setof{v\in[\baire]^N}{(\exists\;t\in\baire) (\exists\;n_0\in\omega)(\forall\;n>n_0)\;(v\rest n,t\rest n)\in T}.$$
\end{df}

\subsection{The rank of a model}

We recall the definition and some properties of a model rank from
\cite{Sh522}. Let $M$ be a model with vocabulary $\tau(M)$ and with the
universe $|M|$. Denote by $M^*$ the collection of all finite sets
$w\subs|M|$ which are {\em independent} in the following sense: if
$a\in w$ and $b_0,\dots,b_{k-1}\in w\setminus\{a\}$ then there is no
quantifier free formula $\phi(x_0,\dots,x_{k-1},y)$ such that
$M\models\phi(b_0,\dots,b_{k-1},a)$ and there exists only finitely many $y\in M$ with $M\models\phi(b_0,\dots,b_{k-1},y)$.

Now define inductively for $w\in M^*$ the {\em rank of} $w$ with respect
to $M$, $\rk(w,M)$, by letting:
\begin{enumerate}
\item[(a)] $\rk(w,M)$ is an ordinal or $\infty$ (where $\al<\infty$
for all ordinals $\al$);
\item[(b)] if $\beta$ is a limit ordinal then $\rk(w,M)\goe\beta$ iff
$\rk(w,M)\goe\al$ for every $\al<\beta$;
\item[(c)] $\rk(w,M)\goe\al+1$ iff for every $a\in w$ and for every
quantifier free formula $\phi\in\tau(M)$ such that
$M\models\phi(a,b_0,\dots,b_{k-1})$, where
$w\setminus\{a\}=\{b_0,\dots,b_{k-1}\}$, there exists $w'\in M^*$ such
that $\rk(w',M)\goe\al$, $w'=\{b_0,\dots,b_{k-1},a,a'\}$, where $a\ne
a'$, and $M\models\phi(a',b_0,\dots,b_{k-1})$.
\end{enumerate} Finally we set
$\rk(w,M)=\al$ iff $\rk(w,M)\goe\al$ and $\rk(w,M)\not\goe\al+1$.
Define $\rk(M)=0$, if $M^*=\emptyset$ and
$$\rk(M)=\sup\{\rk(w,M)+1\colon w\in M^*\},$$ otherwise. Let
$\lam_\al(\aleph_0)$ be the minimal cardinal $\lam$ such that every model
$M$ with countable vocabulary and with universe of size $\goe\lam$, has
rank $\goe\al$. It has been proved in \cite{Sh522} that  $\lam_\al(\aleph_0)\loe\beth_{\al+\omega}$, so in particular $\lam_\al(\aleph_0)$ is well defined.

It is easy to observe that the rank defined above is absolute, but $\lam_\al(\aleph_0)$ may increase in a ZFC
extension of the universe, since new models may be added. It is known (see
\cite{Sh522}) that the sequence
$\{\lam_{\al}(\aleph_0)\}_{\al\loe\omega_1}$ is strictly increasing and
continuous, hence $\aleph_\al\loe\lam_{\al}(\aleph_0)$. It
has been proved in \cite{Sh522} that if $\poset$ is a ccc forcing then in
the generic extension $V^{\poset}$ we have
$\biglam\loe\beth^V_{\omega_1}$. In particular, it is consistent that
$\biglam<2^{\aleph_0}$. We will use the fact that for every $\al<\omega_1$
there exists a model with universe $\lam_\al(\aleph_0)$, with countable
vocabulary and with rank $\al$.

\section{Analytic colorings and homogeneous sets}\label{krowka}

In this section we present some basic properties of colorings and statements of known results.
We are interested in possible cardinalities of homogeneous sets for analytic colorings.
The following lemma allows us to fix the attention to colorings of the Baire space (or the Cantor set).

\begin{lm}\label{kormoran} Let $X$ be a metric space and let $C\subs[X]^N$ be an analytic coloring. Then there exists an analytic coloring $D\subs[\baire]^N$ and a continuous map $\map f{\baire}X$ such that $C=\setof{\img fv}{v\in D}$. If moreover, $X$ is analytic then we may choose $f$ so that $\img f{\baire}=X$. \end{lm}

\begin{pf} By the above remarks, $\suppt(C)$ is analytic. Fix a continuous map $\map f{\baire}X$ such that $\img f{\baire}=\suppt(C)$ (or $\img f{\baire}=X$, if $X$ is analytic). Define $\map\phi{[\baire]^N\times\baire}{[X]^{\loe N}\times\baire}$ by setting $\phi(v,t)=(\img fv,t)$. Then $\phi$ is continuous and we may set $$D=\setof{v\in[\baire]^N}{(\exists\;t)\;(v,t)\in F},$$
where $F=\setof{(v,t)}{\img fv = g(t)}$ and $\map g{\baire}C$ is a fixed continuous map onto $C$. Observe that $F$ is closed, which implies that $D$ is analytic. Clearly $C=\setof{\img fv}{v\in D}$.
\end{pf}

By the above lemma, we will consider coloring trees, which are countable objects defining analytic colorings. This allows us to consider absoluteness of properties of colorings, meaning properties of coloring trees. In particular, observe that the existence of a perfect homogeneous set is absolute, because it can be stated as an absolute property of a coloring tree -- this is also a consequence of the results of  Section \ref{rangakoloru}, where we define an absolute rank function on coloring trees which measures possible cardinalities of homogeneous sets.

The following fact was mentioned, in the context of $\sig$-compact colorings, in \cite{Sh522}.

\begin{prop} The existence of an uncountable homogeneous set for an analytic coloring is absolute between transitive models of ZFC. \end{prop}

\begin{pf} Let $\jezyk(Q)$ be the language $\jezyk$ with additional quantifier $Q$ which means ``there exists uncountably many". We will use a result of Keisler \cite{Kei} which says that $\Theta \subs \jezyk(Q)$ is consistent iff it has a model.
Let $T$ be an $N$-coloring tree of infinite height. There exists $\Theta_0\subs \jezyk$ which describes $T$ (using $\jezyk$ we can represent standard Peano arithmetic, elements of $\baire$, etc.). Let $c$ be an $N$-ary relation which represents that an $N$-tuple belongs to the coloring induced by $T$. Let $R$ be a binary predicate and let $\phi_0(R)$ be ``$R$ is a well-ordering with $\dom(R)\subs \baire$ and for every $x_0,\dots,x_{N-1}\in\dom(R)$ we have $(x_0,\dots,x_{N-1})\in c$". Finally, let $\phi_1(R)$ be ``$(Q\;x)\;x\in \dom(R)$". Let $\Theta=\Theta_0\cup\dn{\phi_0(R)}{\phi_1(R)}$. Then $\Theta$ has a model iff there is an uncountable $C$-homogeneous set, where $C$ is the coloring induced by $T$. By Keisler's theorem, this is equivalent to the fact that $\Theta$ is consistent; the last property is absolute.
\end{pf}

Let us finally mention the following known facts about colorings:

\begin{prop}[\cite{Sh522}, \cite{Kubis1}] Let $C$ be a $G_\delta$-coloring of a Polish space. If there is an uncountable $C$-homogeneous set then there exists also a perfect one. \end{prop}

In \cite{Kubis1} it is proved that if a $G_\delta$-coloring of a Polish space does not contain a perfect homogeneous set then there exists a countable ordinal which bounds the Cantor-Bendixson ranks of all homogeneous sets for this coloring.

\begin{prop}[\cite{Sh522}] Let $\kappa=2^{\aleph_0}$ and let $\poset$ be a Cohen (random) forcing adding $\lam>\kappa$ Cohen (random) reals. Then, in $V^{\poset}$, every analytic coloring which has a homogeneous set of size $\kappa^+$, also has a perfect one. \end{prop}

Thus, it is consistent with ZFC that $2^{\aleph_0}$ is arbitrarily large and every analytic coloring which has a homogeneous set of size $\aleph_2$ also has a perfect one.

It has been proved in \cite{Sh522} that every analytic coloring with a homogeneous set of size $\biglam$ has also a perfect homogeneous set. Moreover, it is consistent that $2^{\aleph_0}$ is arbitrarily large while at the same time $\biglam=\aleph_{\omega_1}$ (recall that $\biglam\goe\aleph_{\omega_1}$ and so far it is not known whether the statement ``$\biglam>\aleph_{\omega_1}$" is consistent). 
In \cite{Sh522}, ccc forcing notions are constructed, which introduce $\sig$-compact pair colorings of the Cantor set with arbitrarily large below $\biglam$, but not perfect, homogeneous sets. In Sections \ref{rangakoloru} and \ref{uniwers} we will improve these results, defining a suitable rank of a coloring and constructing ``universal" colorings of any prescribed rank $\gamma<\omega_1$, which consistently contain homogeneous sets of size $\lam_{\gamma}(\aleph_0)$ but which cannot contain homogeneous sets of size $\lam_{\gamma+1}(\aleph_0)$.

\section{The rank of a coloring}\label{rangakoloru}

In this section we define a rank on coloring trees which measures the ``failures" of constructing perfect homogeneous sets. In a special case of $F_\sig$ colorings of the Cantor set, this rank has already appeared in \cite{Sh522}.

Let $T$ be an $N$-coloring tree. An {\em $n$-approximation} of $T$ is a pair $(v,h)$
such that $v$ is a finite, at least $N$-element subset of $\omega^n$ and 
$h\colon[v]^N\to\omega^n$ is such that $(\vec{v},h(\vec{v}))\in T$ for every $\vec{v}\in[v]^N$. 
Denote by $\app_n(T)$ the set of all $n$-approximations of $T$ and write $\app(T)=\bigcup_{\ntr}\app_n(T)$. For $x\in[\omega^n]^{<\omega}$, $y\in[\omega^m]^{<\omega}$ write $x\lhd y$ if $n<m$ and $x=\setof{t\rest n}{t\in y}$.
Define a strict partial order on $\app(T)$ by setting $(v,h)<(v',h')$ iff $v\lhd v'$ and $h(u)\subs h(u')$ whenever $u\in [v]^N$, $u'\in[v']^N$ and $u\lhd u'$. 

We define the rank $\rk_T(v,h)$ of each $(v,h)\in\app(T)$ as an ordinal or
$\infty$ by letting:
\begin{enumerate}
\item[(i)] $\rk_T(v,h)=\al$ iff $\rk_T(v,h)\goe\al$ and
$\rk_T(v,h)\not\goe\al+1$; if $\beta>0$ is a limit then
$\rk_T(v,h)\goe\beta$ iff $\rk_T(v,h)\goe\al$ for all $\al<\beta$;
\item[(ii)] $\rk_T(v,h)\goe\al+1$ iff for each $p\in v$ there exists
$(v',h')\in\app(T)$ such that $\rk_T(v',h')\goe\al$, $(v,h)<(v',h')$ and $p$ has at
least two extensions in $v'$.
\end{enumerate}
Finally, define $\rk(T)=\sup\{\rk_T(v,h)+1\colon(v,h)\in\app(T)\}$.

Below is proved the crucial property of $\rk_T$.

\begin{prop}\label{ptak} Let $T$ be a coloring tree of infinite height and let $C$ be the coloring induced by $T$. 
If for some $(v,h)\in\app(T)$ we have
$\rk_T(v,h)\goe\omega_1$ then in fact $\rk_T(v,h)=\infty$ and there exists a
perfect $C$-homogeneous set. Conversely, if there exists a perfect $C$-homogeneous set then $\rk(T)=\infty$.
\end{prop}

\begin{pf} First observe that we can find $(v',h')$ satisfying (ii) with
$\rk_T(v',h')\goe\omega_1$ (there are only countably many $(v',h')$ and
uncountably many $\al<\omega_1$). Let $$\beta=\min\{\rk_T(w,g)\colon
(w,g)\in\app(T)\}\setminus\omega_1$$ and fix $(w,g)$ with
$\rk_T(w,g)=\beta$. By the previous argument, for each $p\in w$ there is
$(w',g')>(w,g)$ satifsying (ii) and with $\rk_T(w',g')\goe\omega_1$. It
follows that $\rk_T(w,g)\goe\beta+1$, a contradiction.

Construct a sequence $S=\{(v_n,h_n)\colon\ntr\}$ such
that $(v_n,h_n)<(v_{n+1},h_{n+1})$, each $p\in
v_n$ has at least two extensions in $v_{n+1}$ and $\rk_T(v_n,h_n)=\infty$. To find $(v_{n+1},h_{n+1})$
apply condition (ii) finitely many times, with $\al=\infty$. Define
$$P=\{x\in\baire\colon(\forall\;\ntr)\;x\rest \length(v_n,h_n)\in v_n\}.$$
Then $P$ is a perfect set, by the definition of $S$. It remains to check
that $P$ is $C$-homogeneous. Fix $\vec{x}=\{x_0,\dots,x_{N-1}\}\in[P]^N$.
Let $n_0$ be such that $x_i\rest m$ are pairwise distinct, where
$m=\length(v_{n_0},h_{n_0})$. Let $\vec{v}_n=\{x_i\rest \length(v_n,h_n)\colon i<N\}$.
Then $\vec{v}_n\in[v_n]^N$ for $n\goe n_0$ and $\{h_n(\vec{v}_n)\}_{n\goe
n_0}$ is an increasing chain of initial segments of some $t\in\baire$.
Then $(\vec{x}\rest n,t\rest n)\in T$ for $n\goe m$ and therefore $\vec{x}\in C$. 

For the last statement, assume $P$ is perfect and $C$-homogeneous. Then for every $x\in[P]^N$ there is $\tau(x)\in\baire$ such that $(x,\tau(x))\in[T]$. We say that $(u,h)\in\app(T)$ is {\em induced by} $P$ if there exists $X\in[P]^{<\omega}$ such that $u=\setof{\sig\rest n}{\sig\in X}$ and $h(v)=\tau(y)\rest n$ where $y\subs X$ is such that $v=\setof{\sig\rest n}{\sig\in y}$. Observe that if $(u,h)$ is induced by $P$ then, since $P$ is perfect, there exists $(u',h')>(u,h)$ which is also induced by $P$ and every element of $u$ splits in $u'$.
Thus, by induction on $\beta<\omega_1$ we can prove that $\rk_T(u,h)\goe\beta$ whenever $(u,h)$ is induced by $P$.
\end{pf}

Next we show a connection between our rank and cardinals
$\lam_\gamma(\aleph_0)$.

\begin{prop}\label{liczby} Let $C$ be an analytic coloring of $\baire$ induced by a coloring tree $T$ and let $\gamma\loe\omega_1$. If there exists a
$C$-homogeneous set of size $\lam_\gamma(\aleph_0)$ then
$\rk(C)\goe\gamma$. \end{prop}

\begin{pf} Let $\lam=\lam_\gamma(\aleph_0)$ and let
$\{\eta_\al\colon\al<\lam\}$ be a one-to-one enumeration of a
$C$-homogeneous set of size $\lam$. 
Let $$B=\setof{(v,t)\in[\baire]^N\times\baire}{(\exists\;n_0\in\omega)(\forall\;n>n_0)\; (v\rest n, t\rest n)\in T}.$$
For each
$\{\al_0,\dots,\al_{N-1}\}\in[\lam]^N$ pick $t=t(\al_0,\dots,\al_{N-1})$
such that $(\{\eta_{\al_0},\dots,\eta_{\al_{N-1}}\},t)\in B$.
Now fix $\ntr$, $\vec{v}\in[\omega^n]^N$, $s\in\omega^n$ and define
$$R_{\vec{v},s}=\{(\al_0,\dots,\al_{N-1})\in\lam^N\colon
\vec{v}\lhd\{\eta_{\al_0},\dots,\eta_{\al_{N-1}}\}\Land s\subs
t(\al_0,\dots,\al_{N-1})\}.$$
Now $M=(\lam,\{R_{\vec{v},s}\colon\ntr,\;s\in\omega^n,\;
\vec{v}\in[\omega^n]^N\})$ is a model with countable vocabulary. By the
definition of $\lam$, $\rk(M)\goe\gamma$.

Fix $w\in M^*$, $|w|\goe N$. Let $m$ be such that $v=\{\eta_\al\rest
m\colon\al\in w\}$ has the same size as $w$. For $\vec{v}\in[v]^N$ let
$h(\vec{v})=t(\al_0,\dots,\al_{N-1})\rest m$, where
$\vec{v}=\{\eta_{\al_0}\rest m,\dots,\eta_{\al_{N-1}}\rest m\}$. Then
$(v,h)\in\app_m(C)$; we will say that $(v,h)$ is {\em determined by} $w$.
Now it suffices to prove that $\rk_T(v,h)\goe\rk(w,M)$ whenever $(v,h)$ is
determined by $w$.

Fix $\beta<\gamma$ and assume that this is true whenever $\rk(w,M)<\beta$.
If $\beta$ is a limit then we can apply the continuity of the rank; thus
we may assume that $\beta=\xi+1$. Let $\rk(w,M)\goe\xi+1$ and let
$(v,h)$ be determined by $w$.

For each $\vec{v}\in[v]^N$, $s=h(\vec{v})$, there is
$w(\vec{v})=(\al_0,\dots,\al_{N-1})\in w^N$ ($\al_i$'s are pairwise
distinct) such that $M\models R_{\vec{v},s}(w(\vec{v}))$. Fix $p\in v$. We
have to find $(v',h')$ satisfying condition (ii) of the definition of
$\rk_T$. Let $p=\eta_{\beta}\rest \length(v,h)$, where $\beta\in w$. By the
definition of $\rk(w,M)$, there exists $w'\in M^*$ such that
$\rk(w',M)\goe\xi$, $w'=(w\setminus\{\beta\})\cup\{\beta'\}$ (where
$\beta'\notin w$) and for each $\vec{v}\in[v]^N$ we have $M\models
R_{\vec{v},s}(w'(\vec{v}))$, where $w'(\vec{v})$ is like $w(\vec{v})$,
with $\beta'$ instead of $\beta$. Find $m>\length(v,h)$ and an $m$-approximation
$(v',h')$ determined by $w\cup\{\beta'\}$. Then $(v,h)<(v',h')$ and $p$
has two extensions in $v'$: $\eta_{\beta}\rest m$ and $\eta_{\beta'}\rest
m$. By induction hypothesis, $\rk_T(v',h')\goe\rk(w',M)\goe\xi$. Hence
$\rk(v,h)\goe\xi+1$. \end{pf}

Summarizing the above results we get:

\begin{tw} Let $C$ be an analytic coloring induced by a coloring tree $T$. Then there is no perfect $C$-homogeneous set if and only if $\rk(T)<\omega_1$. If $\gamma<\omega_1$ and there exists a $C$-homogeneous set of size $\lam_\gamma(\aleph_0)$ then $\rk(T)\goe\gamma$. \end{tw}

\begin{wn} Let $C$ be an analytic coloring. If\/ $C$ contains a homogeneous set of size $\lam_\gamma(\aleph_0)$ for every $\gamma<\omega_1$ then $C$ also contains a perfect homogeneous set. \end{wn}

\begin{pf} By Lemma \ref{kormoran}, we can assume that $C$ is a coloring of the Baire space. Now, if $T$ is a coloring tree inducing $C$, then by the above result we have $\rk(T)\goe\omega_1$ and hence there exists a perfect $C$-homogeneous set.
\end{pf}

\begin{uwgi} As we have already mentioned, a $G_\delta$ coloring of a Polish space with an uncountable homogeneous set has also a perfect one. In this case, the rank defined above decides only whether there exists a perfect homogeneous set. For example, if $T$ is a tree inducing a closed coloring $C\subs[\baire]^N$ such that $T$ consists of pairs $(v,0)$ where $v\in[\omega^n]^{<\omega}$ and $0$ is the zero sequence in $\omega^n$, then $\rk(T)$ is not smaller than the supremum of Cantor-Bendixson ranks of $C$-homogeneous sets.
Indeed, let $K$ be $C$-homogeneous with Cantor-Bendixson rank $\goe\gamma$ and denote by $K^{(\al)}$ the $\al$-th derivative of $K$. We claim that for every approximation $(u,h)\in\app(T)$ such that $u\lhd K^{(\al)}$, we have $\rk_T(u,h)\goe\al$. Use induction on $\al$: the statement is true if $\al=0$; if $u\lhd K^{(\al+1)}$ then there exists $u'$ such that $u\lhd u'$, every element of $u$ splits in $u'$ and $u'\lhd K^{(\al)}$. Thus, by induction hypothesis, $\rk_T(u',h')\goe\al$ (where $h'$ is defined in the obvious way), so $\rk_T(u,h)\goe\al+1$. The case of limit $\al$ is the same.
\end{uwgi}

\section{Universal colorings}\label{uniwers}

In this section we prove that for every countable ordinal $\gamma$ there exists an $F_\sig$ $2$-coloring $C$ on the Cantor set, which has rank $\gamma$ and which contains a copy of any $\sig$-compact $2$-coloring of the Cantor set, whose rank is $\loe \gamma$. A similar statement is true for $N$-colorings, where $N>2$, however we will concentrate on $\sig$-compact pair colorings.

\begin{df} A {\em basic coloring tree} is a coloring tree $T$ with the following properties:
\begin{enumerate}
	\item[(i)] $T$ is a pair coloring tree of height $\loe\omega$ (for our purposes we need sometimes finite height).	
$T$ consists of triples $(x,y,k)$ where $x,y\in\omega^n$ for some $\ntr$ and $k\in\omega$ (we write $(x,y,k)$ instead of $(\dn xy, k)$ so $(x,y,k)\in T$ implies $(y,x,k)\in T$). According to the definition of a coloring tree, we denote by $\Lev_n(T)$ the set of all $(x,y,k)\in T$ such that $x,y\in \omega^n$.
	\item[(ii)] If $(x,y,k)\in\Lev_n(T)$ and $m>n$ then there exist $x'\sups x$, $y'\sups y$ such that $(x',y',k)\in \Lev_m(T)$.
	\item[(iii)] For every $n$ the set $\setof{k}{(\exists\;x,y)\;(x,y,k)\in \Lev_n(T)}$ is finite. 
	\item[(iv)] The {\em support tree of\/} $T$, $\suppt(T)=\setof{x}{(\exists\;(x',y,k)\in T)\;x\subs x'}$, determines a compact set (i.e. is finitely branching).
\end{enumerate}
The order of $T$ is defined, as usual, by $(x,y,k)\loe (x',y',k')$ if $x\subs x'$, $y\subs y'$ and $k=k'$. The set $\setof{k}{(\exists\;x,y)\;(x,y,k)\in T}$ will be called the {\em set of colors} of $T$ and denoted by $\colors(T)$. We define $\app(T)$ in the same way like in the general case of coloring trees.
\end{df}

Note that, by condition (iii), we may have that $(x,y,k)\in \Lev_n(T)$ but not necessarily $(x\rest m, y\rest m, k)\in T$ for $m<n$. 
We need condition (iii) only for technical reasons.

\begin{claim} Let $T$ be a basic coloring tree. Then $T$ inudces an $F_\sig$ coloring
$$C=\setof{\dn xy \in[K]^2}{(\exists\;k,n_0)(\forall\;n>n_0)\; (x\rest n,y\rest n, k)\in T},$$ where $K=\suppt(T)\subs \baire$ is compact. Conversely, if $K\subs\baire$ is compact and $C\subs[K]^2$ is $F_\sig$ then $C$ is induced by some basic coloring tree. \end{claim}

\begin{pf} The first statement is clear. For the second one, assume $C=\bigcup_{\ntr} C_n$ is a $2$-coloring of a compact set $K\subs\baire$ such that each $C_n$ is closed. Define $$T=\setof{(x\rest n, y\rest n, k)}{\dn xy \in C_k\Land k\loe n}.$$
Then $T$ is a basic coloring tree which induces $C$. \end{pf}

\begin{df}\label{sarna} Let $\gamma$ be a countable ordinal. A {\em $\gamma$-ranked coloring tree} is a triple $(T,r,c)$, where $T$ is a basic coloring tree, $r$ is a function defined on approximations of $T$ whose values are ordinals $<\gamma$, $c$ is a function defined on $\app(T)$ such that $c(u,h)\in u$ and the following conditions are satisfied:
\begin{enumerate}
	\item[($\mho_1$)] If $(u,h)<(u',h')$ then $r(u,h)\goe r(u',h')$; if moreover $c(u,h)$ splits in $(u',h')$ then $r(u,h)>r(u',h')$.
	\item[($\mho_2$)] If $(u,h)<(u',h')$, $r(u,h)=r(u',h')$ and $|u|=|u'|$ (i.e. $u$ has no splitting in $u'$) then $c(u',h')$ extends $c(u,h)$.
	\item[($\mho_3$)] If $w\subs u$ and $g=h\rest[w]^2$ then $r(w,g)\goe r(u,h)$.
\end{enumerate}
If\/ $\Tau=(T,r,c)$ is a $\gamma$-ranked coloring tree and $n<\Ht(T)$ then we denote by $\Tau\rest n$ the $\gamma$-ranked coloring subtree consisting of $\bigcup_{k<n}\Lev_k(T)$ with the restricted functions $r$ and $c$.
\end{df}

In the above definition, $r$ pretends to be the rank of the coloring induced by
$T$ and $c$ specifies the critical elements (witnesses for the rank).
In fact, $r$ gives an upper bound for the rank.

\begin{claim} Let $(T,r,c)$ be a $\gamma$-ranked coloring tree. Then $\rk_T(u,h)\loe r(u,h)$ for every $(u,h)\in\app(T)$. Conversely, if\/ $T$ is a basic coloring tree of rank $\loe\gamma$ then there exist $r, c$ such that $(T,r,c)$ is a $\gamma$-ranked tree. \end{claim}

\begin{pf} Induction on $r(u,h)$. Assume that $r(u,h)=\beta$ and
$\rk_T(w,g)\loe r(w,g)$ whenever $r(w,g)<\beta$. Suppose
$\rk_T(u,h)>\beta$. Then there exists $(u',h')$ such that $(u,h)<(u',h')$,
$c(u,h)$ splits in $(u',h')$ and $\rk_T(u',h')\goe\beta$. Applying (a)
with induction hypothesis we get $\rk_T(u',h')\loe r(u',h')<\beta$; a
contradiction. For the second part, fix a basic coloring tree $T$ with $\rk(T)\loe\gamma$. Set $r=\rk_T$. Inductively, define $c(u,h)$ to be a critical element in $(u,h)$, taking care of condition (b). 
\end{pf}

\begin{df} (a) Let $\gamma$ be a countable ordinal. A {\em $\gamma$-template} is a $\gamma$-ranked basic coloring tree of height $2$ whose support is a finite binary tree. More specifically, a $\gamma$-template is a triple $(S,r,c)$ such that $S$ is a coloring tree of height $2$, $\suppt(S)\subs\omega^2$ is a finite binary tree and $r$, $c$ are like in the definition of $\gamma$-ranked trees. So $S=S_0\cup S_1$, where $S_i=\Lev_i(S)$, $S$ is finite, each $(x,y,k)\in S_0$ has an extension in $S_1$ and each $x\in\suppt(S_0)$ has at most two extensions in $\suppt(S_1)$. 

 (b) Let $(S,r_S,c_S)$ and $(T,r_T,c_T)$ be $\gamma$-ranked coloring trees. An {\em embedding} of $(S,r_S,c_S)$ into $(T,r_T,c_T)$ is a one-to-one order preserving map $\map f{\suppt(S)}{\suppt(T)}$ such that 
\begin{enumerate}
\item If $x,y$ are on the same level then so are $f(x),f(y)$.
\item $f$ determines a map $\map{f^*}{\colors(S)}{\colors(T)}$ which is one-to-one and such that $(x,y,k)\in S\implies (f(x),f(y),f^*(k))\in T$.
\item For each $(u,h)\in\app(S)$, $f$ determines $(u^f,h^f)\in \app(T)$, given by the formulae $u^f=\img fu$, $h^f(f(x),f(y))=f^*(h(x,y))$. Furthermore, we have $r_S(u,h)\loe r_T(u^f,h^f)$ and $f(c_S(u,h))=c_T(u^f,h^f)$.
\end{enumerate}

 (c) A $\gamma$-ranked coloring tree $(T,r,c)$ is {\em universal} if its height is $\omega$ and for every $\gamma$-template $\Es$, any embedding of $\Es\rest1$ (which will be called a {\em partial embedding of} $\Es$) into $(T,r,c)$ can be extended to a full embedding of $\Es$.
\end{df}

\begin{lm}\label{klucz} Let $\gamma$ be a countable ordinal. There exists a universal $\gamma$-ranked coloring tree. \end{lm}

\begin{pf} Let $\gamma=\bigcup_{\ntr}\Gamma_n$ where $\Gamma_0\subs\Gamma_1\subs\dots$ and each $\Gamma_n$ is finite. We will construct inductively a $\gamma$-ranked coloring tree $\Tau=(T,r,c)$ such that:
\begin{enumerate}
\item[(a)] If $(x,y,k)\in T$ then $(x\concat0, y\concat0,k)\in T$; if $(u,h)\in \app(T)$ then $r(u,h)=r(u\concat0,h')$ and $c(u\concat0,h')=c(u,h)\concat0$, where $h'(x\concat0,y\concat0)=h(x,y)$ for $\dn xy\in[u]^2$.

\item[(b)] If $\Es=(S,r_S,c_S)$ is a $\gamma$-template such that $|\Lev_0(\suppt(S))|< n$, $|\colors(S)|< n$ and $\rng(r_S)\subs \Gamma_n$, then the $0$-th level of $\Es$ can be embedded into $\Tau\rest n$ and any such embedding can be extended to a full embedding into $\Tau\rest (n+1)$.
\end{enumerate}
Suppose $\Tau\rest n$ has been defined. We need to extend $T$ to the $n$-th level, preserving the above conditions. First, add all elements of the form $x\concat0$, where $x\in \Lev_{n-1}(\suppt(T))$ and extend $r$ and $c$ according to (a). Formally, we extend $T$ by adding $T_0=\setof{(x\concat0,y\concat0,k)}{(x,y,k)\in \Lev_{n-1}(T)}$.

Enumerate as $f_1,\dots,f_{l}$ the set of all (isomorphism types of) partial embeddings of simple $\gamma$-templates into the so far defined $T$. Let $(S_i,r_i,c_i)$ denote the $\gamma$-template corresponding to $f_i$, i.e. $\dom(f_i)=\Lev_0(\suppt(S_i))$. Assume that the range of $f_i$ lies on the $n$-th level of $T$ -- otherwise, by condition (a) we can lift it to the $n$-th level (adding zeros) without changing anything. Without loss of generality we may assume that
$$\setof{s\concat0}{s\in\dom(f_i)}\subs \Lev_1(\suppt(S_i))\subs \setof{s\concat\eps}{s\in\dom(f_i)\Land \eps\in\dn0i}.$$ 
Note that together with $f_i$ we have also a map $f_i^*$ which maps colors of the $0$-th level of $S_i$ to $\colors(T)$. Extend $f_i^*$ to a one-to-one map $\map{f_i^{**}}{\colors(S_i)}{W}$ where $W\sups \colors(T)$ is such that colors coming from the $1$-st level are mapped to $W\setminus\colors(T)$. Now define
$$T_i=\setof{(f_i(x)\concat\eps, f_i(y)\concat\delta,f_i^{**}(k))}{(x\concat\eps,y\concat\delta,k)\in \Lev_1(S_i)}.$$
Let $T'=T\cup\bigcup_{i<l+1}T_i$. Observe that $(u,h)\in\app_n(T')$ iff $(u,h)\in\app_n(T\cup T_i)$ for some $i\loe l$, because a pair of the form $(x\concat\eps,y\concat\delta)$, where $x,y\in\Lev_{n-1}(\suppt(T))$, has a color in $T'$ only if $\eps,\delta\in \dn0i$ for some $i$. 

Extend each $f_i$ to $f_i'$ naturally (i.e. setting $f_i'(x\concat\eps)=f_i(x)\concat\eps$) and define $r(u,h)$ and $c(u,h)$ in such a way that $f_i'$ becomes an embedding. It is clear that conditions ($\mho_1$)--($\mho_3$) are satisfied.

Finally, we need to further extend $T$ in order to fulfill the first part of condition (b) for $n+1$, i.e. every $\gamma$-template $\Es=(S,r_S,c_S)$ with $|\Lev_0(\suppt(S))|<n+1$, $|\colors(S)|<n+1$ and with $\rng(f_S)\subs\Gamma_{n+1}$ can be partially embedded into the $n$-th level of $(T,r,c)$.
There exist only finitely many isomorphism types of such $\gamma$-templates. Let $\Es_0,\dots,\Es_{k-1}$ enumerate all these types. For each $i<k$ choose a finite set $K_i\subs \omega$ in such a way that $K_i\cap K_j=\emptyset$ whenever $i<j$ and $|K_i|=|\Lev_0(\suppt(\Es_i))|$.  

Without loss of generality, we may assume that $\Lev_0(\suppt(\Es_i))=K_i$ and $\colors(\Es_i)\cap\colors(T)=\emptyset$. Denote by $\theta$ the constant-zero sequence in the $(n-1)$-th level of $\Tau\rest n$ and define $g_i(s)=\theta\concat s$ and extend $T$ by adding all triples of the form $(\theta\concat s,\theta\concat t,k)$, where $(s,t,k)\in \Lev_0(\Es_i)$ (and extend $r$ $c$ to the new approximations according to $g_i$). Now each $g_i$ becomes a partial embedding of $\Es_i$. 
It is clear that conditions ($\mho_1$) -- ($\mho_3$) are satisfied. This finishes the proof.
\end{pf}

\begin{lm}\label{krab} Let $(T,r,c)$ be a universal $\gamma$-ranked tree. Then every $\gamma$-ranked coloring tree can be embedded into $(T,r,c)$. 
\end{lm}

\begin{pf} Let $\Es=(S,r_S,c_S)$ be a $\gamma$-ranked coloring tree. We can assume without loss of generality that $\suppt(S)$ is a binary tree: embed $\suppt(S)$ into a binary tree $B$ and then, using this embedding, define a new $\gamma$-ranked coloring tree whose support is $B$. So, assuming $\suppt(S)$ is binary, for every $\ntr$ the tree $\Lev_n(\Es)\cup\Lev_{n+1}(\Es)$ together with restricted $r_S$ and $c_S$ is a $\gamma$-template. Using the universality of $(T,r,c)$, embed the $0$-th level of $\Es$ into $(T,r,c)$.

Now, assume that $f$ is a partial embedding of $\Es$ into $(T,r,c)$ defined on the first $n$ levels of $\Es$. Denote by $\Es'$ the coloring tree consisting of $\Lev_n(\Es)\cup\Lev_{n+1}(\Es)$. Then $\Es'$ is a $\gamma$-template (because $\Es$ is binary) and $f$ is a partial embedding of $\Es'$. By universality, $f$ can be extended to a full embedding of $\Es'$ and this extension defines an embedding of $\bigcup_{i\loe n+1}\Lev_i(\Es)$. 
Thus, by induction we can construct a full embedding of $\Es$ into $(T,r,c)$.
\end{pf}

Now we are able to state the main results of this section.

\begin{tw} Let $C$ be an $F_\sig$ pair coloring of the Cantor set $2^\omega$ which is induced by a universal $\gamma$-ranked coloring tree. Then for every coloring $D\subs[2^\omega]^2$ induced by a basic coloring tree of rank $\loe\gamma$ there exists a topological embedding $\map\phi{2^\omega}{2^\omega}$ such that $\dn{\phi(x)}{\phi(y)}\in C$ for every $\dn xy\in D$. \end{tw}

\begin{pf} Let $(T,r,c)$ denote the universal $\gamma$-ranked tree which induces $C$. Then $\suppt(T)$ is a finitely branching subtree of $\omega^{<\omega}$. Let $Y=[\suppt(T)]$ (the set of branches through $\suppt(T)$), so $Y$ is a compact subspace of $\baire$ homeomorphic to a subset of the Cantor set. Moreover $C\subs [Y]^2$. 
Let $S$ be a basic coloring tree of rank $\loe\gamma$ which induces $D$. Let $X$ be a subspace of $\baire$ which is homeomorphic to the Cantor set and contains $[\suppt(S)]$. Fix $x_0\in X$ and fix a new color $c_0$ (i.e. $c_0\notin \colors(S)$). Define $S'=S\cup\setof{(\dn{x_0\rest n}{x\rest n},c_0)}{\ntr,\;x\in X\setminus\sn{x_0}}$. Then $\suppt(S')=X$ and $S'$ induces an $F_\sig$ coloring containing $D$.
Furthermore, $\rk(S')=\rk(S)$ because if $(u,h)$, $(u',h')$ are in $\app(S')$ and $\dn{s}{s_0}\in u$ are such that $s_0\subs x_0$ and $h(s,s_0)=c_0$, then $s_0$ cannot split in $(u',h')$.
Thus, without loss of generality we may assume that $[\suppt(S)]=X$. 
Now by Lemma \ref{krab}, $S$ can be embedded into $T$, which means that there exists a topological embedding $\map \phi XY$ with the property that $\dn{\phi(x)}{\phi(y)}\in C$ whenever $\dn xy\in D$. 
\end{pf}

\begin{tw}\label{main1} Let $\gamma<\omega_1$ and let $C$ be a coloring induced by a universal $\gamma$-ranked coloring tree. Let $\lam\loe\lam_\gamma(\aleph_0)$. Then there exists a ccc forcing notion $\poset$ which forces that $C$ contains a homogeneous set of size $\lam$. Specifically, $|\poset|=\lam$ and there exists a family $\Dee$ of dense subsets of $\poset$ with $|\Dee|=\lam$ and such that a $\Dee$-generic filter defines a $C$-homogeneous set of size $\lam$.
\end{tw}

\begin{pf} Fix a model $M$ with countable vocabulary, with universe $\lam$ 
and with rank $\loe\gamma$. For each $w\in [M]^{<\omega}$ fix a formula 
$\phi=\phi_M(w)$ and an element $x=c_M(w)$ which witness that 
$\rk(w,M)\not\goe\xi+1$, where $\xi=\rk(w,M)$. We will call $\phi_M(w)$ 
and $c_M(w)$ the {\em critical formula} and the {\em critical element} of 
$w$. Fix a universal $\gamma$-ranked coloring tree $(T,r,c)$ inducing $C$. Let $\poset$ consist 
of quadruples of the form $p=(\eta,w,n,g)$,where: \begin{enumerate}

\item[(i)] $w$ is a finite subset of the universe of $M$, $\eta$ is a 
one-to-one function with domain $w$ and with range contained in the $n$-th 
level of $\suppt(T)$; $g$ is a function from $[w]^2$ into the set of colors of
$T$.

\item[(ii)] If $\al<\beta$ are in $w$ then $(\eta(\al), \eta(\beta),
g(\al,\beta)) \in T$.

\item[(iii)] For every $v\subs w$ we have $\rk(v,M)\loe r(\eta[v],g^v)$
where, by definition, $g^v(\eta(\al),$ $\eta(\beta))=g(\al,\beta)$;
moreover $\eta(c_M(v))=c(\eta[v],g^v)$.

\end{enumerate}
We will write $(\eta^p,w^p,n^p,g^p)$ to emphasize $p$. We order $\poset$
by the following rule: $p\loe q$ iff $n^p\loe n^q$, $w^p\subs w^q$, 
$\eta^q(\al)\rest n^p = \eta^p(\al)$ for $\al\in w^p$ and $g^p = 
g^q\rest[w]^2$ (recall that a bigger element is a stronger condition).

\begin{claim}\label{jeden} $\poset$ is ccc. \end{claim}

\begin{pf}
Fix an uncountable set $\Cee\subs\poset$. Shrinking $\Cee$, we may assume
that $n^{p}=n$ for every $p\in\Cee$ and that each two elements $p,q\in
\Cee$ are isomorphic, i.e. there exists an order preserving
bijection $f\colon w^{p}\to w^{q}$ such that
\begin{enumerate}

\item[(a)] $\eta^{q}(f(\al)) = \eta^{p}(\al)$ and
$g^{q}(f(\al),f(\beta)) = g^{p}(\al,\beta)$ for $\al,\beta\in
w^{p}$;

\item[(b)] $c_M(v)=c_M(f[v])$ and $\phi_M(v)=\phi_M(f[v])$ for $v\subs 
w^{p}$.

\end{enumerate}
Let $\Cee_0\subs \Cee$ be uncountable such that $\{w^p\colon p\in
\Cee_0\}$ forms a $\Delta$-system with root $v$. Now we claim that each
two elements of $\Cee_0$ are compatible. Fix $p,q\in \Cee_0$ and let $f$
be an isomorphism as described above. We are going to define a $\gamma$-template $(S,r^S,c^S)$. We start with describing $\suppt(S)$. The
$0$-th level of $\suppt(S)$ will be $w^p$, the $1$-st level will be 
$w^p\concat0\cup(w^p\setminus v)\concat1$. The order of $\suppt(S)$ is 
natural. Then $f$ defines a correspondence between $w^p\cup w^q$ and 
$\Lev_1(\suppt(S))$ given by $x\mapsto x\concat0$ for $x\in w^p$ and $x\mapsto 
f^{-1}(x)\concat1$ for $x\in w^q\setminus v$. Denote this correspondence 
by $\chi\colon w^p\cup w^q\to \Lev_1(\suppt(S))$.
Now define $S_0=\Lev_0(S)$ by 
assigning colors to all pairs of $\Lev_0(\suppt(S))$ according to $g^p$, i.e. 
$(x,y,g^p(x,y))\in S$ (so we have unique 
colors on level $0$). Define 
\begin{align*}
\Lev_1(S)=&\{(x\concat0,y\concat0,k)\colon (x,y,k)\in 
S_0\}\cup\{(x\concat0,y\concat1,g^p(x,y))\colon x\in v,\;y\in w^p\setminus 
v\}\\ \cup&\{(x\concat0,y\concat1,k_*)\colon x,y\in w^p\setminus v\},
\end{align*}
where $k_*\notin \rng(g^p)=\rng(g^q)$ is fixed arbitrarily. For 
$(u,h)\in\app_0(S)$ define $r^S(u,h)=\rk(u,M)$ and $c^S(u,h)=c_M(u)$.
Similarly, define $r^S,c^S$ for $(u,h)\in\app_1(S)$ by letting 
$r^S(u,h)=\rk(u',M)$ and $c^S(u,h)=c_M(u')$ where $u'=\chi^{-1}[u]$. 
Observe that $(S,r^S,c^S)$ is a $\gamma$-template: if $(u,h)<(u',h')$ in 
$S$ and $c^S(u,h)$ splits in $(u',h')$, then $u'$ is of the form 
$v_0\concat0\cup\{x\concat0,x\concat1\}$ where $v_0\subs v$ and $x\in 
w^p\setminus v$ (if $x\ne y$ are in $w^p\setminus v$ then the pair 
$x\concat0,y\concat1$ has different color than $x,y$). Now, as 
$c^S(u)=c_M(v_0\cup\{x\})$ and we have
$\phi_M(v_0\cup\{x\})=\phi_M(v_0\cup \{f(x)\})$, we have 
$\rk(v_0\cup\{x,f(x)\})<\rk(v_0\cup\{x\},M)=\rk(v_0\cup\{f(x)\},M)$ by the 
definition of $\phi_M$ and $c_M$. It follows that $r^S(u',h')<r^S(u,h)$.

Now it is clear that $\eta^p,\eta^q$ and $g^p,g^q$ define an embedding
$(S,r^S,c^S)\rest1 \to (T,r,c)$ which, by universality, can be extended to 
an embedding $\psi\colon (S,r^S,c^S)\to (T,r,c)$. This allows us to define a 
condition extending both $p$ and $q$. Indeed, define $w=w^p\cup w^q$, 
define $\eta\colon w\to \suppt(T)$ by setting $\eta(\al)=\psi_0(\chi(\al))$; 
finally define $g(\al,\beta)=\psi_1(\chi(\al),\chi(\beta))$. Then 
$t=(\eta,w,n,g)\in\poset$, where $n$ is such that 
$\psi_0[\Lev_1(\suppt(S))]\subs\Lev_n(\suppt(T))$. Clearly, we have $p,q\loe t$.
\end{pf}

\begin{claim}\label{dwa} For every $\al<\lam$ and $\ntr$ the set
$$D_{\al,n}=\{p\in \poset \colon \al\in w^p \Land n^p\goe n\}$$ is dense in
$\poset$. \end{claim}

\begin{pf}
Fix $p=(\eta^p,w^p,n^p,g^p) \in \poset$. We will define a simple 
$\gamma$-template $(S,r^S,c^S)$ by setting 
$$\suppt(S)=w^p\cup\{x\concat0\colon x\in w^p\}\cup\{x_0\concat 1\},$$ where 
$x_0\in w^p$ is fixed arbitrarily. Next define colors for $\Lev_0(\suppt(S))$ 
according to $g^p$ and take a new color for each pair of the form
$\{x\concat 0, x_0\concat 1\}$. Define $(r^S,c^S)$ using the rank and 
critical elements taken from $M$. Now, $g^p$ defines a partial embedding 
of $(S,r^S,c^S)$ to $(T,r,c)$ which, by universality, can be extended to a 
full embedding $\phi\colon (S,r^S,c^S) \to (T,r,c)$. This allows us to 
define a condition $q$ which is stronger than $p$, $\al\in w^q$ and 
$n^q>n^p$. Repeating this procedure finitely many times, we can get 
$n^q\goe n$. \end{pf}

Now let $G$ be $\poset$-generic. Define $$\eta_\al=\bigcup\{\eta^p(\al)\colon
p\in G\Land \al \in w^p\}.$$ By Claim \ref{dwa} the definition is correct.
Each $\eta_\al$ is an element of the Cantor set and if $\al<\beta$ then
$\eta_\al\ne \eta_\beta$ by (i). Let $g(\al,\beta)=g^p(\al,\beta)$ where $p\in
G$ is such that $\al,\beta\in w^p$. This does not depend on the choice of $p$.
We have $(\eta_\al\rest n, \eta_\beta\rest n, g(\al,\beta)) \in T$ for big
enough $n$. This means that $\{\eta_\al\colon \al<\lam\}$ is a $C$-homogeneous
set of size $\lam$. This completes the proof.
\end{pf}

\begin{wn} Let $C$ be a coloring induced by a universal $\gamma$-ranked tree $(T,r,c)$. Then $\rk(T)=\gamma$ and for every tree $S$ which induces $C$ we have $\rk(S)\goe\gamma$. \end{wn}

\begin{wn} Assume $MA_\kappa$. Then for every $\gamma<\omega_1$, every coloring induced by a universal $\gamma$-ranked tree contains a homogeneous set of cardinality $\min\dn\kappa{\lam_\gamma(\aleph_0)}$. \end{wn}

\section{An application: defectedness of subsets of linear spaces}

Let $X$ be a subset of a Polish linear space and let $N>1$. Define $$C=\setof{S\in[X]^N}{\conv S\not\subs X},$$ where $\conv S$ denotes the convex hull of $S$, i.e. the smallest convex set containing $S$.
Then $C$ is called the {\em $N$-defectedness coloring} of $X$. $C$-homogeneous sets are called {\em $N$-cliques in} $X$. The existence of perfect cliques has been studied, mainly for closed sets, in \cite{KPS, Ko1, FKo, GKKS, GK}. This is related to the study of {\em uncountably convex sets}, i.e. sets which are not countable unions of convex sets. Several examples of closed uncountably convex sets are known, for which the minimal number of convex sets whose union is the given set, is consistently smaller than the continuum. Of course, such sets cannot contain perfect cliques so the induced defectedness colorings have countable rank.

If $X$ is a closed set, then the defectedness is open, so either there exists a perfect clique in $X$ or else all cliques are countable (and the supremum of their Cantor-Bendixson ranks is countable). It is proved in \cite{Kubis1} that an analytic subset of the plane either contains a perfect $3$-clique or for every $N>1$, all $N$-cliques in $X$ are countable. The case of higher dimensions is not known. Below we show that the defectedness of a Borel set is analytic, so the results of the previous sections can be applied here. Next we prove that any $\sig$-compact $N$-coloring of the Cantor set can be realized as defectedness in some $G_\delta$ set in $\Err^{N+1}$. It has been shown in \cite{KPS,GKKS, GK} that certain well known clopen colorings of the Cantor set can be realized as defectedness in some compact subsets of finite-dimensional real linear spaces.

\begin{prop} Let $X$ be a Borel subset of a Polish linear space. Then for every $N>1$, the $N$-defectedness coloring is analytic. \end{prop}

\begin{pf} Fix a continuous map $f\colon\baire\to X$ onto $X$. Identify
$[X]^N$ with $X^N$ and consider $\phi\colon\Sigma_N\times X^N\to E$
(where $E$ is the linear space containing $X$) such that $\Sigma_N$ is the
standard simplex spanned by the canonical linear base in $\Err^N$ and
$\phi$ assigns to $(t,x)$ the convex combination of $x$ with coefficients
in $t$. Define $$C=\{x\in(\baire)^N\colon(\exists\;t)\;|f(x)|=N\Land
\phi(t,f(x))\in E\setminus X\}.$$
Then $C$ is the projection of some Borel subset of
$(\baire)^N\times\Sigma_N$, hence it is analytic. Finally, the $N$-defectedness of $X$ is a continuous image of $C$ (by using $f$).
\end{pf}

\begin{tw} Let $N>1$, let $E$ be a Polish linear space of dimension at
least $2N-1$ and let $C\subs[2^{\omega}]^N$ be an $F_\sig$ coloring. Then there
exists a $G_\delta$ set $S\subs E$ and a homeomorphism
$f\colon2^\omega\to X$, where $X\subs S$, $S\setminus X$ is a countable union of
convex sets and for each $T\in[X]^N$ we have $\conv T\not\subs S$ iff $f^{-1}[T]\in
C$. \end{tw}

\begin{pf} We assume that $E=\Err^{2N-1}$ with the Euclidean metric. For
$x\in[E]^N$ let $I(x)=(\conv x)\setminus\{\conv y\colon y\in[x]^{N-1}\}$.
Then $I(x)\nnempty$ iff $x$ is affinely independent.
Construct a Cantor set
$X=\bigcap_{\ntr}\bigcup_{s\in2^n}D_s$ on the unit sphere of $E$, in such
a way that for each $\ntr$:
\begin{enumerate}
\item[(1)] no $2N$ points taken from distinct $D_s$'s ($s\in2^n$) are
contained in a $(2N-2)$-dimensional affine subspace;
\item[(2)] if $a_0,a_1\in[2^k]^N$ are distinct and $x_0,x_1$ are $N$-element
sets such that $x_i$ is selected from the family $\{D_s\colon s\in a_i\}$,
$i<2$, then $I(x_0)\cap I(x_1)=\emptyset$.
\end{enumerate} This is possible, because the dimension of $E$ is $2N-1$. In fact, (2) follows from (1), because if $I(x_0)\cap I(x_1)\nnempty$ then the affine hull of $x_0\cup x_1$ has dimension at most $(|x_0|-1)+(|x_1|-1)\loe2N-2$.

This construction defines a homeomorphism $f\colon2^\omega\to X$.
Observe that if $x_0,x_1\in[X]^N$ are distinct then $I(x_0)\cap
I(x_1)=\emptyset$, by (2).
Let $C=\bigcup_{\ntr}C_n$ be an $N$-coloring of $2^\omega$, where each
$C_n$ is closed. Fix $m<\omega$ and $x\in[X]^N$. Let $c=\sum_{s\in x}\frac 1Ns$. For each $s\in x$ let $p(s)$ be the point in the segment $[c,s]$ whose distance from $s$ equals $\min\dn{1/m}{\|s-c\|}$. Finally, let $b_m(x)=\setof{p(s)}{s\in x}$. Observe that: 
\begin{enumerate}
	\item[($*$)] if $p_n\in b_m(x_n)$ and $q=\lim_{n\to\infty}p_n$ then $q\in
b_m(x)$, where $x$ is the limit of some convergent subsequence of
$\{x_n\}_{\ntr}$.
\end{enumerate}
Define $S_m=(\conv X)\setminus
F_m$, where $F_m=\{b_m(f[\vec{x}])\colon\vec{x}\in C_m\}$.
Observe that $S_m$ is open in $\conv X$, by ($*$).
Finally, let $S=\bigcap_{m<\omega}S_m$. It is clear that $C$ determines
the $N$-defectedness of $S$ restricted to $X$. It remains to show that
$S\setminus X$ is a countable union of convex sets. If $p\in S\setminus X$ then some
neighborhood of $p$ is disjoint from all $F_m$'s: if not, then $p$ is an
accumulation point of a sequence $\{p_n\}_{\ntr}$ such that $p_n\in F_n$;
but $\dist(p_n,X)\loe1/n$ so $p\in X$, a contradiction. Hence, $p$ has a
convex neighborhood in $S\setminus X$. Now, as $S\setminus X$ has a
countable base, it follows that it consists of countably many convex sets. \end{pf}

\end{document}